\documentclass[12pt,reqno]{amsart}
\usepackage{enumerate}
\usepackage{amsmath, amsfonts, amssymb, amsthm}
\def\N{{\Bbb N}}
\def\Z{{\Bbb Z}}

\def\l{\left}
\def\r{\right}
\def\ls{\leqslant}
\def\gs{\geqslant}

\def\({\bigg(}
\def\[{\bigg[}
\def\){\bigg)}
\def\]{\bigg]}

\def\t{\text}
\def\f{\frac}
\def\mo{{\rm mod}}

\def\eq{\equiv}

\numberwithin{equation}{section}
\newtheorem{Theorem} {Theorem} [section]

\newtheorem{Corollary} [Theorem] {Corollary}
\newtheorem{Conjecture}[Theorem]{Conjecture}
\theoremstyle{definition}

\newtheorem{Remark}[Theorem]{Remark}

\begin{document}
\hbox{In honor of Prof. M. B. Nathanson on the occasion of his 65th
birthday.}
\medskip
\title
[Mixed sums of primes and other terms]{Mixed sums of primes and
other terms}
\author[Zhi-Wei Sun]
{Zhi-Wei Sun}

\thanks{Research supported by the National Natural Science Foundation of China (grant 10871087).}
\address{Department of Mathematics, Nanjing University, Nanjing 210093, People's Republic of China}
\email{zwsun@nju.edu.cn\quad{\it Homepage}:\ {\tt http://math.nju.edu.cn/$\sim$zwsun}}


\begin{abstract} In this paper we study mixed sums of primes and
linear recurrences. We show that if $m\eq2\ (\mo\ 4)$ and $m+1$ is a
prime then $(m^{2^n-1}-1)/(m-1)\not=m^n+p^a$ for any $n=3,4,\ldots$
and prime power $p^a$. We also prove that if $a>1$ is an integer,
$u_0=0$, $u_1=1$ and $u_{i+1}=au_i+u_{i-1}$ for $i=1,2,3,\ldots$,
then all the sums $u_m+au_n\ (m,n=1,2,3,\ldots)$ are distinct. One
of our conjectures states that any integer $n>4$ can be written as
the sum of an odd prime and two positive Fibonacci numbers.
\end{abstract}

\maketitle

\section{Introduction}

Let us first recall the famous Goldbach conjecture in additive number
theory.

\begin{Conjecture}[Goldbach's Conjecture]\label{Goldbach}
Any even integer $n\gs4$ can be written as the sum of two primes.
\end{Conjecture}

The number of primes not exceeding $n\geqslant2$ is approximately
$n/\log n$ by the prime number theorem. Hardy and Littlewood
conjectured that the number of ways to write an even integer $n\gs
4$ as the sum of two primes is given asymptotically by
$$\f {cn}{\log^2n}\prod_{p\mid n}\l(1+\f1{p-2}\r),$$
where $c=2\prod_p(1-(p-1)^{-2})=1.3203\cdots$ is a constant  and $p$
runs over odd primes. (Cf. \cite[pp. 159-164]{Gu}.)

Goldbach's conjecture remains open, and the best result in this
direction is Chen's theorem (cf. \cite{C}): Each large
even integer can be written as the sum of a prime and a product of at most two
primes.

Those integers $T_x=x(x+1)/2$ with $x\in\N=\{0,1,2,\ldots\}$ are called triangular numbers.
There are less than $\sqrt{2n}$ positive triangular numbers below an integer $n\gs2$,
so triangular numbers are more sparse than prime numbers.
In 2008 the author made the following conjecture.

\begin{Conjecture}[Sun \cite{S09}]\label{pT} {\rm (i)} Each natural number $n\not=216$
can be written in the form $p+T_x$ with $x\in\N$, where $p$ is zero or a prime.

{\rm (ii)} Any odd integer greater than $3$ can be written in the form $p+x(x+1)$, where $p$ is a prime and
$x$ is positive integer.
\end{Conjecture}

Douglas McNeil (University of London) (cf. \cite{M2})
 has verified parts (i) and (ii) up to $10^{10}$ and $10^{12}$ respectively.
 The author \cite{S-offer} would like to offer 1000 US dollars for the first positive solutions to both (i) and (ii),
 and \$200 for the first explicit counterexample to (i) or (ii).

Powers of two are even much more sparse than triangular numbers. In a letter to
Goldbach, Euler posed the problem whether any odd integer $n>1$ can
be expressed in the form $p+2^a$, where $p$ is a prime and
$a\in\N$. This question was reformulated by
Polignac in 1849. By introducing covers of the integers by residue
classes, Erd\H os \cite{E50} showed that there exists an infinite
arithmetic progression of positive odd integers no term of which is
of the form $p+2^a$. (See also Nathanson \cite[pp. 204-208]{N}.) On
the basis of the work of Cohen and Selfridge \cite{CS}, the author
\cite{S00} proved that  if
$$x\eq 47867742232066880047611079\ (\mo\ M)$$
with
 $$\aligned M=&2\times3\times5\times7\times11\times13\times17\times19\times31\times37
 \\&\times41\times61\times73\times97
\times109\times151\times241\times257\times331
\\=&66483084961588510124010691590,
\endaligned$$ then $x$ is not of the form
$\pm p^a\pm q^b$ where $p,q$ are primes and $a,b\in\N$.

In 1971 Crocker \cite{Cr} proved that there are infinitely many positive odd integers
not of the form $p+2^a+2^b$ where $p$ is a prime and $a,b\in\Z^+=\{1,2,3,\ldots\}$.
Here are the first few such numbers greater than 5 recently found by Charles Greathouse (USA):
$$6495105,\ 848629545,\ 1117175145,\ 2544265305,\ 3147056235,\ 3366991695.$$
Note that 1117175145 even cannot be written in the form $p+2^a+2^b$ with $p$ a prime and $a,b\in\N$.

Erd\H os (cf. \cite{E97}) asked whether there is a positive integer
$k$ such that any odd number greater than 3 can be written the sum
of an odd prime and at most $k$ positive powers of two. Gallagher
\cite{G} proved that for any $\varepsilon>0$ there is a positive
integer $k=k(\varepsilon)$ such that those positive odd integers not
representable as the sum of a prime and $k$ powers of two form a
subset of $\{1,3,5,\ldots\}$ with lower asymptotic density at least
$1-\varepsilon$. In 1951 Linnik \cite{L} showed that there exists a
positive integer $k$ such that each large even number can be written
as the sum of two primes and $k$ positive powers of two; Heath-Brown
and Puchta \cite{HP} proved that we can take $k=13$. (See also Pintz
and Ruzsa \cite{PR}.)

In March 2005  Georges Zeller-Meier \cite{ZM} asked whether $2^{2^n-1}-2^n-1$ is composite for every $n=3,4,\ldots$.
Clearly an affirmative answer follows from part (i) of our following theorem in the case $m=2$.

\begin{Theorem}\label{p22} {\rm (i)} Let $m\eq 2\ (\mo\ 4)$ be an integer with $m+1$ a prime.
Then, for each $n=3,4,\ldots$, we have
$$\f{m^{2^n-1}-1}{m-1}\not=m^n+p^a,$$
where $p$ is any prime and $a$ is any nonnegative integer.

{\rm (ii)} Let $m$ and $n$ be integers greater than one. Then
$$\f{m^{2^n}-1}{m-1}\not=p+m^a+m^b,$$
where $p$ is any prime, $a,b\in\N$ and $a\not=b$.
\end{Theorem}

\begin{Remark} In the case $m=2$, part (ii) of Theorem \ref{p22} was observed by  A. Schinzel and Crocker independently
in the 1960s, and this plays an important role in Crocker's result about
$p+2^a+2^b$. In 2001 the author and Le \cite {SL} proved that for
$n=4,5,\ldots$ we cannot write $2^{2^n-1}-1$ in the form
$p^\alpha+2^a+2^b$, where $p$ is a prime, $a,b,\alpha\in\N$ and
$a\not=b$.
\end{Remark}

For any integer $m>1$, the sequence $\{m^n\}_{n\gs0}$ is a
first-order linear recurrence with earlier terms dividing all later
terms. To seek for good representations of integers, we'd better
turn resort to second-order linear recurrences whose general term
usually does not divide all later terms.

The famous Fibonacci sequence $\{F_n\}_{n\gs0}$ is defined as follows:
$$F_0=0,\ F_1=1,\ \t{and}\ F_{n+1}=F_n+F_{n-1}\ \t{for}\ n=1,2,3,\ldots.$$
Here are few initial Fibonacci numbers:
$$F_0=0<F_1=F_2=1<F_3=2<F_4=3<F_5=5<F_6=8<F_7=13<F_8=21<\cdots.$$
It is well known that
$$F_n=\f1{\sqrt5}\(\(\f{1+\sqrt5}2\)^n-\(\f{1-\sqrt5}2\)^n\)\ \ \ \t{for all}\ n\in\N.$$
Clearly $F_n<2^{n-1}$ for $n=2,3,\ldots$, and
$$F_n\sim \f{\varphi^n}{\sqrt5}\qquad\ (n\to+\infty),$$
where
$$\varphi=\f{1+\sqrt5}2=1.618\cdots.$$
Note that $2\mid F_n$ if and only if $3\mid n$.

It is not known whether the positive integers not of the form $p+F_n$ with $p$ a prime and $n\in\N$
form a subset of $\Z^+$ with positive lower asymptotic density. However, Wu and Sun \cite{WS} were able to construct
a residue class containing no integers of the form $p^a+F_{3n}/2$ with $p$ a prime and $a,n\in\N$. Note that $u_n=F_{3n}/2$
is just half of an even Fibonacci number; also $u_0=0$, $u_1=1$, and $u_{n+1}=4u_n+u_{n-1}$ for $n=1,2,3,\ldots$.

On December 23, 2008 the author \cite{pFF} formulated the following conjecture.
\begin{Conjecture}[Conjecture on Sums of Primes and Fibonacci Numbers]\label{pFF}
Any integer $n>4$ can be written as the sum of an odd prime and two positive Fibonacci numbers.
We can require further that one of the two Fibonacci numbers is odd.
\end{Conjecture}

\begin{Remark}\label{R-pFF} For a large integer $n$, there are about $\log n/\log\varphi$
Fibonacci numbers below $n$ but there are about $n/\log n$ primes
below $n$. So, Fibonacci numbers are much more sparse than prime
numbers and hence the above conjecture looks more difficult than the
Goldbach conjecture. D. McNeil (cf. \cite{M2, M3}) has
verified Conjecture \ref{pFF} up to $10^{14}$. The author (cf. \cite{S-offer}) would like
to offer 5000 US dollars for the first positive solution published
in a well-known mathematical journal and \$250 for the first
explicit counterexample which can be rechecked by the author via
computer. Note that Conjecture \ref{pFF} implies that for any odd prime $p$ we can find an odd prime $q<p$ such that
$p-q$ can be written as the sum of two odd Fibonacci numbers.
\end{Remark}

Recall that the Pell sequence $\{P_n\}_{n\gs0}$ is defined as follows.
$$P_0=0,\ P_1=1,\ \t{and}\ P_{n+1}=2P_n+P_{n-1}\ \ \t{for}\ n=1,2,3,\ldots.$$
It is well known that
$$P_n=\f1{2\sqrt2}\l((1+\sqrt2)^n-(1-\sqrt2)^n\r)\ \ \ \t{for all}\ n\in\N.$$
Clearly $P_n>2^{n}$ for $n=6,7,\ldots$, and
$$P_n\sim \f{(1+\sqrt2)^n}{2\sqrt2}\qquad\ (n\to+\infty).$$

On Jan. 10, 2009, the author \cite {FF} posed the following conjecture which is an analogue of Conjecture \ref{pFF}.
\begin{Conjecture}[Conjecture on Sums of Primes and Pell Numbers]\label{pP2P}
Any integer $n>5$ can be written as the sum of an odd prime, a Pell number and twice a Pell number.
We can require further that the two Pell numbers are positive.
\end{Conjecture}

\begin{Remark}\label{R-pP2P} D. McNeil (cf. \cite{S-offer}) has verified Conjecture \ref{pP2P} up to $5\times10^{13}$ and found no counterexample.
The author (cf. \cite{S-offer}) would like to offer 1000 US dollars for the first positive solution published in a well-known mathematical journal
and \$100 for the first explicit counterexample
which can be rechecked by the author via computer.
\end{Remark}

 Soon after he learned Conjecture  \ref{pP2P} from the author, Qing-Hu Hou (Nankai University) observed (without proof) that
 all the sums $P_s+2P_t\ (s,t=1,2,3,\ldots)$ are distinct.
 Clearly Hou's observation follows from our following theorem.

\begin{Theorem}\label{uau} Let $a>1$ be an integer, and set
$$u_0=0,\ u_1=1, \ \t{and}\ u_{i+1}=au_i+u_{i-1}\ \t{for}\ i=1,2,3,\ldots.$$
Then no integer $x$ can be written as $u_m+au_n$ (with $m\in\N$ and $n\in\Z^+$) in at least two ways,
except in the case $a=2$ and $x=u_0+au_2=u_2+au_1=4$.
\end{Theorem}

\begin{Remark}\label{R-uau} Note that if $n\in\Z^+$ then $u_{n+1}+au_0=au_n+u_{n-1}$.
\end{Remark}

\begin{Corollary}\label{P+2P} Let $k,l,m,n\in\Z^+$. Then $P_k+2P_l=P_m+2P_n$ if and only if $k=m$ and $l=n$.
\end{Corollary}

\begin{Remark}\label{R-P+2P} In view of Corollary \ref{P+2P},
we can assign an ordered pair $\langle m,n\rangle\in\Z^+\times\Z^+$
the code $P_m+2P_n$. Recall that a sequence $a_1<a_2<a_3<\cdots$ of
positive integers is called a Sidon sequence if all the sums of
pairs, $a_i+a_j$, are all distinct. An unsolved problem of Erd\H os
(cf. \cite[p. 403]{Gu}) asks for a polynomial $P(x)\in\Z[x]$ such
that all the sums $P(m)+P(n)\ (0\ls m<n)$ are distinct.
\end{Remark}

Motivated by Conjecture \ref{pFF} and its variants, Qing-Hu Hou and Jiang Zeng (University of Lyon-I)
formulated the following conjecture during their visit to the author in Jan. 2009.

\begin{Conjecture}[Hou and Zeng \cite{HZ}]\label{pFC}
Any integer $n>4$ can be written as the sum of an odd prime, a positive Fibonacci number and a Catalan number.
\end{Conjecture}

\begin{Remark}\label{R-pP2P} Catalan numbers are integers of the form
$$C_n=\f1{n+1}{2n\choose n}={2n\choose n}-{2n\choose n+1}\quad (n\in\N),$$
which play important roles in combinatorics (see, e.g., Stanley
\cite[Chapter 6]{St}). They are also determined by $C_0=1$ and the
recurrence
$$C_{n+1}=\sum_{k=0}^nC_kC_{n-k}\quad(n=0,1,2,\ldots).$$
By Stirling's formula, $C_n\sim 4^n/(n\sqrt{n\pi})$ as $n\to+\infty$.
D. McNeil \cite{M3} has verified Conjecture \ref{pFC} up to $3\times10^{13}$ and found no counterexample.
Hou and Zeng  would like to offer 1000 US dollars for the first positive solution published in a well-known mathematical journal
and \$200 for the first explicit counterexample
which can be rechecked by them via computer. Note that 3627586 cannot be written in the form $p+2F_s+C_t$ with $p$ a prime and $s,t\in\N$.
\end{Remark}

The Lucas sequence $\{L_n\}_{n\gs0}$ is defined as follows.
$$L_0=2,\ L_1=1,\ \t{and}\ L_{n+1}=L_n+L_{n-1}\ (n=1,2,3,\ldots).$$
It is known that
$$L_n=2F_{n+1}-F_n=\(\f{1+\sqrt5}2\)^n+\(\f{1-\sqrt5}2\)^n$$
for every $n=0,1,2,3,\ldots$.

On Jan. 16, 2009 the author (cf. \cite{FF}) made the following conjecture which is
similar to Conjecture \ref{pFC}.

\begin{Conjecture}\label{pLC}
Each integer $n>4$ can be written as the sum of an odd prime, a
Lucas number and a Catalan number.
\end{Conjecture}
\begin{Remark}\label{R-pLC} D. McNeil \cite{M3} has verified Conjecture \ref{pLC}
up to $10^{13}$ and found no counterexample. Note that 1389082 cannot be written
in the form $p+2L_s+C_t$ with $p$ a prime and $s,t\in\N$.
\end{Remark}

Recall that there are infinitely many positive odd integers not of the form $p+2^a+2^b$ with $p$ a prime and $a,b\in\Z^+$.
However, Crocker's trick in his proof of this result does not work for the form $p+2^a+k2^b$ with $p$ a prime and $a,b\in\Z^+$,
where $k$ is an odd integer greater than one.
On Jan. 21, 2009 the author (cf. \cite{FF}) made the following conjecture.

\begin{Conjecture}[Conjecture on Sums of Primes and Powers of Two]\label{p222}
Any odd integer greater than $8$ can be written as the sum of an odd prime and three positive powers of two.
Moreover, we can write any odd integer $n>10$ in the form $p+2^a+3\times 2^b=p+2^a+2^b+2^{b+1}$ with $p$ a prime and $a,b\in\Z^+$.
\end{Conjecture}

\begin{Remark}\label{p2k2} The author verified Conjecture \ref{p222} for odd integers below $10^7$. Later, on the request of the author,
Qing-Hu Hou and Charles Greathouse  continued the verification for odd integers below $2\times 10^8$ and $10^{10}$ respectively.
 Note that if $k>61$ is odd then $2k+127$ cannot be written in the form $p+2^a+k2^b$
with $p$ an odd prime and $a,b\in\Z^+$ since $3+2+k2^2>2k+127$ and $127$ is not of the form $p+2^a$.
For $k\in\{3,5,\ldots,61\}\setminus\{47,51\}$, the author (cf. \cite{pFF}) checked odd integers below $10^8$
and found no odd integer $n>2k+3$ not of the form
$p+2^a+k2^b$ with $p$ an odd prime and $a,b\in\Z^+$.
\end{Remark}

 We are going to prove Theorems \ref{p22} and \ref{uau} in the next section.
Section 3 is devoted to our discussion of Conjecture \ref{pFF} and its variants.

\section{Proofs of Theorems \ref{p22} and \ref{uau}}

\medskip
\noindent{\it Proof of Theorem \ref{p22}}.
For $n=2,3,\ldots$ we clearly have
$$\aligned(m-1)\prod_{k=0}^{n-1}\l(m^{2^k}+1\r)=&\l(m^{2^0}-1\r)\l(m^{2^0}+1\r)\l(m^{2^1}+1\r)\cdots\l(m^{2^{n-1}}+1\r)
\\=&\l(m^{2^1}-1\r)\l(m^{2^1}+1\r)\cdots\l(m^{2^{n-1}}+1\r)
\\=&\cdots=\l(m^{2^{n-1}}-1\r)\l(m^{2^{n-1}}+1\r)=m^{2^n}-1.
\endaligned$$

(i) Fix an integer $n\gs 3$. Write $n+1=2^kq$ with $k\in\N$, $q\in\Z^+$ and $2\nmid q$. Since
$$2^n=(1+1)^n\gs 1+n+\f{n(n-1)}2>n+1,$$
we must have $0\ls k\ls n-1$. Thus $m^{2^k}+1$ divides both $(m^{2^n}-1)/(m-1)$ and
$m^{n+1}+1=(m^{2^k})^q+1$.
Set
$$d_n=\f{m^{2^n-1}-1}{m-1}-m^n.$$
Then
$$md_n=\f{m^{2^n}-m}{m-1}-m^{n+1}=\f{m^{2^n}-1}{m-1}-(m^{n+1}+1)$$
and hence $m^{2^k}+1$ divides $d_n$.

Suppose that $d_n$ is a prime power. By the above, we can write $d_n=p^a$, where $a\in\N$ and
$p$ is a prime divisor of  $m^{2^k}+1$. As $m$ is even, $p$ is an odd prime.
Since
$$m^{p-1}\eq1\ (\mo\ p)\ \ \t{and}\ \ m^{2^{k+1}}\eq(-1)^2=1\ (\mo\ p),$$
we have
$$m^{\gcd(p-1,2^{k+1})}\eq1\ (\mo\ p).$$
But
$$m^{2^k}\eq-1\not\eq1\ (\mo\ p),$$
so $p\eq1\ (\mo\ 2^{k+1})$.
Note that
$$p^a=\f{m^{2^n-1}-1}{m-1}-m^n=\sum_{k=0}^{2^n-2}m^k-m^n\eq 1+m+m^2\ (\mo\ m^3).$$
If $k>0$, then $p\eq1\ (\mo\ 2^2)$ and hence
$$p^a\eq1\not\eq 1+m\ (\mo\ 2^2),$$
which contradicts the congruence $p^a\eq1+m\ (\mo\ m^2)$.
So $k=0$, $p\mid m^{2^0}+1$ and hence $p=m+1$. (Recall that $m+1$ is a prime.)
It follows that $p^a$ is congruent to $1$ or $m+1$ modulo 8. Since $1+m+m^2\not\eq1, m+1\ (\mo\ 8)$,
we get a contradiction.
This proves part (i).

(ii) Let $a>b\gs0$ be integers with $m^a+m^b<(m^{2^n}-1)/(m-1)$. Clearly $2^n>a>b$. Write $a-b=2^kq$ with $k\in\N$, $q\in\Z^+$ and $2\nmid q$.
Then $0\ls k<n$ and hence $d=m^{2^k}+1$ divides both $(m^{2^n}-1)/(m-1)$ and $m^{a-b}+1=(m^{2^k})^q+1$. Thus
$$\f{m^{2^n}-1}{m-1}-m^a-m^b$$
is a multiple of $d$. Observe that
$$\aligned \f{m^{2^n}-1}{m-1}=&\f{m^{2^n-2}-1}{m-1}\l(m^{2^n-2}+1\r)\l(m^{2^n-1}+1\r)
\\>&\l(m^{2^n-2}+1\r)\l(m^{2^n-1}+1\r)\gs(m^b+1)(m^{a-b}+1)\gs m^a+m^b+d.
\endaligned$$
So $d$ is a proper divisor of $D=(m^{2^n}-1)/(m-1)-m^a-m^b$. This shows that $D$ cannot be a prime.
We are done. \qed

\medskip
\noindent {\it Proof of Theorem \ref{uau}}. Observe that
$$u_0=0<u_1=1<u_2=a<u_3<u_4<\cdots.$$
By induction,
$$u_{2i}\eq u_0=0\ (\mo\ a)\ \t{and}\ u_{2i+1}\eq u_1=1\ (\mo\ a)\quad\t{for}\ i=0,1,2,\ldots.$$
We will make use of these simple properties.

Let $k,m\in\N$ and $l,n\in\Z^+$ with $k\ls m$. Below we discuss the
equation $u_k+au_l=u_m+au_n$.

{\it Case} 1. $k=m$.

 In this case,
 $$u_k+au_l=u_m+au_n\Rightarrow u_l=u_n\Rightarrow l=n.$$

{\it Case} 2. $k=l<m$.

If $k=l<m-1$ then
$$u_k+au_l<u_{m-2}+au_{m-1}=u_m<u_m+au_n.$$
When $k=l=m-1$, as $u_m\not\eq u_{m-1}\ (\mo\ a)$ we have
$$u_k+au_l=(a+1)u_{m-1}\not=u_m+au_n.$$

{\it Case} 3. $l<k<m$.

In this case,
$$u_k+au_l\ls u_k+au_{k-1}<au_k+u_{k-1}=u_{k+1}\ls u_m<u_m+au_n.$$

{\it Case} 4. $k<l<m$.

In this case,
$$u_k+au_l\ls au_l+u_{l-1}=u_{l+1}\ls u_m<u_m+au_n.$$

{\it Case} 5. $k<m\ls l$.

Suppose that $u_k+au_l=u_m+au_n$. Then
$$u_l>\f{u_k+au_l-u_m}a=u_n\gs u_l-\f{u_m}a\gs\f{a-1}au_l\gs (a-1)u_{l-1}\gs u_{l-1}.$$
It follows that
$$k=0,\ m=l=2, \ \t{and}\ u_n=u_{l-1}=u_1=1.$$
Thus $au_2=u_2+au_1$, i.e., $a^2=a+a$ and hence $a=2$.

Combining the above we have completed the proof. \qed

\begin{Remark}\label{F+F=F+F} By modifying the proof of Theorem \ref{uau}, we can determine all the solutions
of the equation $F_k+F_l=F_m+F_n$ with $k,l,m,n\in\N$.
\end{Remark}

\section{Discussion of Conjecture \ref{pFF} and its variants}

Concerning Conjecture \ref{pFF}, we mention that there are very few natural numbers not representable
as the sum of a prime $p\eq5\ (\mo\ 6)$ and two Fibonacci numbers.
Bjorn Poonen (MIT) informed the author that by a heuristic argument there should be infinitely
many positive integers not in the form $p+F_s+F_t$ if we require that the
prime $p$ lies in a fixed residue classe with modulus greater than one.
McNeil \cite{M1, M3} made a computer search to find natural numbers not representable as
the sum of a prime $p\eq5\ (\mo\ 6)$, an odd Fibonacci number and a positive Fibonacci number; he found that there
are totally 729 such numbers in the interval $[0,10^{14}]$, 277 of which (such as 857530546) even cannot be written
as the sum of a prime $p\eq5\ (\mo\ 6)$ and two Fibonacci numbers.

In 2008 the author (cf. \cite{pFF, FF}) also made the following conjecture which is similar to Conjecture \ref{pFF}.
\begin{Conjecture} \label{pFFk}
{\rm (i)} Any integer $n>4$ can be written as the sum of an odd prime, a positive Fibonacci number
and the square of a positive Fibonacci number. We can require further that one of the two Fibonacci numbers is odd.

{\rm (ii)} Each integer $n>4$ can be written as the sum of an odd prime, a positive Fibonacci number
and the cube of a positive Fibonacci number. We can require further that one of the two Fibonacci numbers is odd.
\end{Conjecture}

\begin{Remark}\label{R-pFFk} Note that 900068 cannot be written as the sum of a prime, a Fibonacci number and the fourth power
of a Fibonacci number.
Also,
$$F_n^3\sim\f{\varphi^{3n}}{(\sqrt5)^3}=\f{(4.236\cdots)^n}{5\sqrt5}\ \quad(n\to +\infty).$$
\end{Remark}

Let $k\in\{1,2,3\}$. For $n\in\Z^+$ let $r_k(n)$ denote the number of ways to write $n$ as the sum of an odd prime, a positive Fibonacci number
and the $k$th power of a positive Fibonacci number with one of the two Fibonacci numbers odd. That is,
$$r_k(n)=|\{\langle p,s,t\rangle:\, p+F_s+F_t^k=n,
\ p\ \t{is an odd prime},\ s,t\ge2, \ \t{and}\ 2\nmid F_s\ \t{or}\ 2\nmid F_t\}|.$$
The author has investigated values of the quotient
$$s_k(n)=\f{r_k(n)}{\log n}$$
via computer, and conjectured  that
$$c_k=\liminf_{n\to+\infty}s_k(n)>0.$$
Numerical data suggest that $2<c_1<3$.
In fact, the author computed all values of $s_1(n)$ with $10^{50}\ls n\ls 10^{50}+4\times10^4$, and here are the two smallest values:
$$s_1(10^{50}+39030)=2.22359\cdots\ \ \t{and}\ \ s_1(10^{50}+5864)=2.29037\cdots.$$

Here is another variant of Conjecture \ref{pFF} made by the author (cf. \cite{pFF,FL}).
\begin{Conjecture} \label{pFL}
{\rm (i)} Any integer $n>4$ can be written as the sum of an odd prime, an odd Lucas number and a positive Lucas number.
For $k=2,3$ we can write any integer $n>4$ in the form $p+L_s+L_t^k$, where $p$ is an odd prime, $s,t\gs0$,
and $L_s$ or $L_t$ is odd.

{\rm (ii)} Each integer $n>4$ can be written as the sum of an odd prime, a positive Fibonacci number
and twice a positive Fibonacci number (or half of a positive Fibonacci number). We can also represent
any integer $n>4$ as the sum of an odd prime, twice a positive Fibonacci number,
and the square of a positive Fibonacci number.

{\rm (iii)} Any integer $n>4$ can be written in the form $p+F_s+L_t$
with $p$ an odd prime, $s>0$, and $F_s$ or $L_t$ odd.
\end{Conjecture}

\begin{Remark}\label{R-pFL} The author verified Conjectures \ref{pFFk} and \ref{pFL} for $n\ls 3\times10^7$.
Qing-Hu Hou found that 17540144 cannot be written as the sum of a prime, a Lucas number and the fourth power
of a Lucas number. McNeil (cf. \cite{M2}) has verified the first assertions
in parts (i) and (ii) of Conjectures  \ref{pFFk} and \ref{pFL} up to $10^{12}$. He (cf. \cite{M3}) has also verified
part (iii) of Conjecture \ref{pFL} up to $10^{13}$, and found that 36930553345551 cannot be written as the sum of a prime,
a Fibonacci number and an even Lucas number.
\end{Remark}

What about the representations $n=p+P_s+kP_t$ with $k\in\{1,3,4\}$ related to
Conjecture \ref{pP2P}? Note that 2176 cannot be written as the sum of a prime and two Pell numbers.
McNeil \cite{M3} found that
393185153350 cannot be written as the sum of a prime, a Pell number and three times a Pell number,
and the smallest integer greater than 7 not representable as the sum of a prime, a Pell number
and four times a Pell number is
$$872377759846\approx 8.7\times 10^{11}.$$
The companion Pell sequence $\{Q_n\}_{n\gs0}$ is defined by
$$Q_0=Q_1=2\ \t{and}\ Q_{n+1}=2Q_n+Q_{n-1}\ (n=1,2,3,\ldots).$$
McNeil \cite{M3} found that the smallest integer greater than 5 not representable as the sum of a prime,
a Pell number and a companion
Pell number is 169421772576.

McNeil's counterexamples seem to suggest that Conjecture \ref{pP2P} might also have large counterexamples.
However, in the author's opinion, the large counterexamples to the representations $n=p+P_s+3P_t$
and $n=p+P_s+4P_t$ hint that they are very close to the ``truth" (Conjecture  \ref{pP2P}).
Corollary \ref{P+2P} is also a good evidence to support Conjecture  \ref{pP2P}.
To expel suspicion, the author has investigated the behavior of the representation function
$$r(n)=|\{\langle p,s,t\rangle:\ p+P_s+2P_t=n \ \t{with}\ p\ \t{a prime}\ \t{and}\ s,t\gs0\}|.$$
For $n\in[10^{50}, 10^{50}+10081]$ most values of $s(n)=r(n)/\log n$ lies in the interval $(1,2)$,
the smallest value of $s(n)$ with $n$ in the range is
$$s(10^{50}+10045)=\frac{76}{\log(10^{50}+10045)}\approx 0.66.$$
The author also computed the values of $s(n)$ with $n\in[10^{200},10^{200}+100]$,
the smallest value and the largest value are
$$s(10^{200}+33)=\frac{443}{\log(10^{200}+33)}\approx 0.96$$
and
$$s(10^{200}+18)=\frac{824}{\log(10^{200}+18)}\approx 1.79$$
respectively. The author conjectured that
$$c=\liminf_{n\to+\infty}s(n)\in(0.6,1.2).$$

\bigskip

\noindent{\bf Acknowledgment}. The author wishes to thank Dr. Douglas McNeil who has checked almost all conjectures
mentioned in this paper (on the author's request) via his quite efficient and powerful computation.


\bigskip

\begin{thebibliography}{20}


\bibitem{C} Jing-run Chen, {\it On the representation of a large even integer as the sum of a prime and the product
of at most two primes},
Sci. Sinica {\bf 16} (1973), 157--176.

\bibitem{CS} F. Cohen and J. L. Selfridge, {\it Not every number is the sum or difference of two prime powers},
Math. Comp. {\bf 29} (1975), 79--81.

\bibitem{Cr} R. Crocker, {\it On a sum of a prime and two powers of two},
Pacific J. Math. {\bf 36} (1971), 103--107.

\bibitem{E50} P. Erd\H os, {\it On integers of the form $2^k+p$ and some related problems},
Summa Brasil. Math. {\bf 2} (1950), 113--123.

\bibitem{E97} P. Erd\H os, {\it Some of my favorite problems and results},
in: The Mathematics of Paul Erd\H os, I (R. L. Graham and J. Ne\v set\v ril, eds.), Algorithms and Combinatorics 13,
Springer, Berlin, 1997, pp. 47--67.

\bibitem{G} P. X. Gallagher, {\it Primes and powers of 2}, Invent. Math. {\bf 29} (1975), 125--142.

\bibitem{Gu} R. K. Guy, {\it Unsolved Problems in Number Theory}, 3rd Edition, Springer, New York, 2004.

\bibitem{HP} D. R. Heath-Brown and J.-C. Puchta, {\it Integers represented as a sum of primes and powers of two},
Asian J. Math. {\bf 6} (2002), 535--565.

\bibitem{HZ} Q. H. Hou and J. Zeng, Sequences A154404 in On-Line Encyclopedia of Integer Sequences,
\newline{\tt http://www.research.att.com/$\sim$njas/sequences/A154404}

\bibitem{L} Yu. V. Linnik, {\it Prime numbers and powers of two}, Trudy Mat. Inst. Steklov.
{\bf 38} (1951), 152--169.

\bibitem{M1} D. McNeil, {\it Sun's strong conjecture} (a message to Number Theory Mailing List in Dec. 2008),
\newline {\tt http://listserv.nodak.edu/cgi-bin/wa.exe?A2=ind0812\&L=nmbrthry\&T=0\&P=3020}.

\bibitem{M2} D. McNeil, {\it Various and sundry} (a message to Number Theory Mailing List in Jan. 2009),
\newline {\tt http://listserv.nodak.edu/cgi-bin/wa.exe?A2=ind0901\&L=nmbrthry\&T=0\&P=840}.

\bibitem{M3} D. McNeil, Private communications in Jan. 2009.

\bibitem{N} M. B. Nathanson, {\it Additive Number Theory: The Classical Bases}, Grad. Texts in Math., Vol. 164,
Springer, New York, 1996.

\bibitem{PR} J. Pintz and I. Z. Ruzsa,  {\it On Linnik's appproximation to Goldbach's problem, I}, Acta Arith.
{\bf 109}(2003), 169--194.

\bibitem{St} R. P. Stanley, {\it Enumerative Combinatorics},
Vol. 2, Cambridge Univ. Press, Cambridge, 1999.

\bibitem{S00}  Z. W. Sun, {\it On integers not of the form $\pm p^a\pm q^b$},
Proc. Amer. Math. Soc. {\bf 128} (2000), 997--1002.

\bibitem{S09}  Z. W. Sun, {\it On sums of primes and triangular numbers},
Journal of Combinatorics and Number Theory {\bf 1} (2009), 65--76.

\bibitem{pFF}  Z. W. Sun, Four messages to the Number Theory Mailing List,
 \newline {\tt
http://listserv.nodak.edu/cgi-bin/wa.exe?A2=ind0812\&L=nmbrthry\&T=0\&P=2140}
\newline {\tt
http://listserv.nodak.edu/cgi-bin/wa.exe?A2=ind0812\&L=nmbrthry\&T=0\&P=2704}
\newline {\tt
http://listserv.nodak.edu/cgi-bin/wa.exe?A2=ind0812\&L=nmbrthry\&T=0\&P=3124}
\newline {\tt
http://listserv.nodak.edu/cgi-bin/wa.exe?A2=ind0901\&L=nmbrthry\&T=0\&P=1886}


\bibitem{FF}  Z. W. Sun, Sequences A154257, A154258, A154263, A154536, A154940, A155860 in On-Line Encyclopedia of Integer Sequences,
\newline {\tt http://www.research.att.com/$\sim$njas/sequences/A154257}
\newline {\tt http://www.research.att.com/$\sim$njas/sequences/A154258}
\newline {\tt http://www.research.att.com/$\sim$njas/sequences/A154263}
\newline {\tt http://www.research.att.com/$\sim$njas/sequences/A154536}
\newline {\tt http://www.research.att.com/$\sim$njas/sequences/A154940}
\newline {\tt http://www.research.att.com/$\sim$njas/sequences/A155860}

\bibitem{FL}  Z. W. Sun, Sequences A154285, A154290, A154417, A155114 in On-Line Encyclopedia of Integer Sequences,
\newline{\tt http://www.research.att.com/$\sim$njas/sequences/A154285}
\newline {\tt http://www.research.att.com/$\sim$njas/sequences/A154290}
\newline {\tt http://www.research.att.com/$\sim$njas/sequences/A154417}
\newline {\tt http://www.research.att.com/$\sim$njas/sequences/A155114}

\bibitem{S09}  Z. W. Sun, {\it On sums of primes and triangular numbers},
Journal of Combinatorics and Number Theory {\bf 1} (2009), 65--76.

\bibitem{S-offer}  Z. W. Sun, {\it Offer prizes for solutions to my main conjectures involving primes}
(a message to the Number Theory Mailing List in Jan. 2009),
\newline {\tt
http://listserv.nodak.edu/cgi-bin/wa.exe?A2=ind0901\&L=nmbrthry\&T=0\&P=1395}


\bibitem{SL} Z. W. Sun and M. H. Le, {\it Integers not of the form $c(2^a+2^b)+p^\alpha$}, Acta Arith.
{\bf 99} (2001), 183--190.

\bibitem{WS} K. J. Wu and Z. W. Sun, {\it Covers of the integers with odd moduli and their applications to the forms
$x^m-2^n$ and $x^2-F_{3n}/2$},
Math. Comp., in press. {\tt arXiv:math.NT/0702382}

\bibitem{ZM} G. Zeller-Meier, {\it Not prime for each $n\gs 3$} (a message to Number Theory Mailing List in March 2005),
{\tt http://listserv.nodak.edu/cgi-bin/wa.exe?A2=ind0503\&L=nmbrthry\&T=0\&P=2173}.


\end{thebibliography}
\end{document}